\documentclass[11pt]{article}
\usepackage{amsmath,amssymb}

\usepackage{epsfig}
\usepackage{hyperref}
\ifx\onlineVersion\undefined%
\long\def\online#1#2{#2}%
\else%
\long\def\online#1#2{#1}%
\fi%
\newdimen\appletheight
\newdimen\appletwidth
\def\applet#1#2#3{%
\appletwidth=#2%
\appletheight=#3%
\vbox to\appletheight{
\vfill\hbox to\appletwidth{%
\vbox to \appletheight{}%
\hfill}}%
}

\newtheorem{theorem}{Theorem}[section]
\newtheorem{proposition}[theorem]{Proposition}
\newtheorem{lemma}[theorem]{Lemma}
\newtheorem{definition}[theorem]{Definition}
\newtheorem{corollary}[theorem]{Corollary}

\begin{document}
%
\begin{center}
{\large\bf
Conformally symmetric circle packings.\\
 A generalization of Doyle spirals}
\end{center}
\vspace{0.5cm}

\begin{center}
{\sc Alexander I.\,Bobenko}\footnote{E--mail:
\href{mailto:bobenko@math.tu-berlin.de}{{\tt bobenko} {\makeatother @ \makeatletter}{\tt math.tu-berlin.de }}} and
{\sc Tim \,Hoffmann}\footnote{E--mail:
\href{mailto:timh@sfb288.math.tu-berlin.de}{{\tt timh} {\makeatother @ \makeatletter}{\tt sfb288.math.tu-berlin.de }}}
\end{center}
\begin{center}
Fachbereich Mathematik, Technische Universit\"at Berlin, \\
Stra{\ss}e des 17 Juni 136, 10623 Berlin, Germany
\end{center}
\vspace{0.5cm}



\section{Introduction}
Circle packings (and more generally patterns) as discrete analogs of conformal
mappings is a fast developing field of research on the border of analysis and
geometry.  Recent progress was initiated by Thurston's idea \cite{T} about
the approximation of the Riemann mapping by circle packings. The corresponding
convergence was proven by Rodin and Sullivan \cite{RS}; many additional connections with analytic
functions, such as the discrete maximum principle and Schwarz's lemma \cite{R}, the discrete
uniformization theorem \cite{BS}, etc., have emerged since then.

The topic ``circle packings'' is also a natural one for computer experimentation and visualization.
Computer experiments demonstrate a surprisingly close analogy of the
classical theory in the emerging ``discrete analytic
function theory'' \cite{DS}.
Although computer experiments give convincing evidence for the existence of
discrete analogs of many standard holomorphic functions, the Doyle spirals (which are discrete analogs
of the exponential function, see section \ref{s.Doyle}) are the only circle packings described
explicitly.

Circle packings are usually described analytically in the Euclidean setting, i.e. through their radii
function. On the other hand, circles and
the tangencies are preserved by the fractional-linear transformations of the Riemann sphere
(M\"obius transformations). It is natural to study circle packings in this setting, i.e. modulo the
group of the M\"obius transformations.
He and Schramm \cite{HS} developed a conformal description of hexagonal circle packings, which helped
them to show that Thurston's convergence of hexagonal circle packings to the Riemann mapping is
actually $C^\infty$. They describe circle packings in terms of the cross-ratios
$$
q(a,b,c,d):= \frac{(a-b)(c-d)}{(b-c)(d-a)}
$$
of their touching points.
\begin{figure}[t]
  \begin{center}
    \online{
       \centerline{\applet{cp.FlowerappletLite.class}{280bp}{280bp}}
     \caption{A circle flower. You can drag the magenta colored
       points. The blue one on the center circle is determined by the multi-ratio
       condition. Right mouse click over the circles gives you more
       viewer options. For help on the viewer click \protect\hyperlink{applethelp}{here}.}%
    }{%
    \input aflower1.tex
    \caption{A circle flower.}
    }\label{f.flower}%
  \end{center}
\end{figure}

In \cite{S} Schramm introduced circle patterns with the combinatorics of the square grid (SG patterns).
In many aspects the SG theory is analogous to the theory of the hexagonal circle packings. However,
the SG theory is analytically simpler.
The corresponding discrete equations describing the SG patterns, in the Euclidean as well as in the
conformal setting, turn out to be integrable \cite{BP}. Methods of the theory of integrable equations
allowed us in \footnote{Discrete
$z^2$ and $\log z$ have been conjectured by Schramm and Kenyon earlier \cite{WWW}.} \cite{AB} to find
Schramm's circle patterns which are analogs of the holomorphic functions $z^\alpha, \log z$.

One big question is which results on the Schramm's circle patterns carry over to the hexagonal setting,
in particular
whether some discrete standard functions can be described explicitly. This is closely related to
the question of integrability of the basic discrete equations for hexagonal circle packings (He-Schramm
equation, see section \ref{s.Analytic}).
In the present paper the first simple step in this direction is made.
We study (surprisingly non-trivial) conformal geometry of hexagonal circle packings. In terms of this
approach, a  special class of {\em conformally symmetric} circle packings, which are
generalizations of the Doyle spirals, is introduced and all such packings are described explicitly.

\online{}{Since this paper deals with families of circle packings it
  seems natural to show not only arbitrarily choosen members in the
  figures, but to provide a possibility to present them all. Therefore
  there is an interactive version of this paper available \cite{online}. It has some
  of the figures replaced by applets, that allow to explore the
  families directly. See section~\ref{s.online} for more information
  on this version.
  }{}
\vspace{2mm}

\noindent {\bf Acknowledgments.}
The authors thank U. Hertrich-Jeromin, U. Pinkall, Yu. Suris and E. Tjaden for helpful discussions.
\section{Geometry of circle flowers and conformally symmetric circle packings}  \label{s.Geometry}

This paper concerns patterns of circles in the plane called hexagonal circle
packings. Their basic unit is the {\em flower}, consisting of a {\em center}
circle tangent to and surrounded by {\em petals}. A {\em hexagonal} flower is
illustrated in Fig.\ref{f.flower}; the six petals form a closed chain which wraps once in
the positive direction about the center. Whereas the neighboring petals touch,
the circles of not-neighboring petals of a flower may intersect. 
We call a flower
{\em immersed} if none of its
circles degenerates to a point. A hexagonal circle packing is a collection of
oriented circles where each of its internal circles is the center of a hexagonal flower.
Orientations of the circles should agree: at the touching points the orientations
of the touching circles must be opposite.
A hexagonal circle packing can be labeled by the triangular (hexagonal) lattice
$$
HL=n+m e^{i\pi/3}\in {\mathbb C}, \qquad n,m\in {\mathbb Z}.
$$
or by one of its subset. A circle packing is called {\em immersed} if all
its flowers are immersed. Immersions of the whole $HL$ are called {\em entire}.
Fractional-linear transformations of the complex plane (M\"obius transformations)
preserve circles, their orientation and their tangencies. In this paper we study circle
packings modulo the group of M\"obius transformations.

The center circle of a flower contains 6 points $z_1,\ldots, z_6$ (see
Fig.\ref{f.flower}) where it
touches the petals. We call them {\em center touching points} of a circle flower.

\begin{proposition}                 \label{p.points-geometry}
Let $z_1,\ldots, z_6$ be cyclicly ordered\footnote{Positive orientation of the ordering is
assumed.} points on a circle $C$. Then the
following three statements are equivalent:
\begin{itemize}
\item [(i)] There exists a flower with the center $C$ and center touching points
$z_1,\ldots, z_6$,
\item [(ii)] The multi-ratio $m$ of $z_1,\ldots, z_6$ is equal to $-1$, i.e.
\begin{equation}                            \label{m-ratio}
m(z_1,z_2,z_3,z_4,z_5,z_6):=
\frac{(z_1-z_2)(z_3-z_4)(z_5-z_6)}{(z_2-z_3)(z_4-z_5)(z_6-z_1)}=-1,
\end{equation}
\item [(iii)] There exists an involutive M\"obius transformation $M$
  (M\"obius involution) such that
$$
M (z_k)=z_{k+3} \qquad\qquad (k\ {\rm mod}\ 6).
$$
\end{itemize}
\end{proposition}
\begin{figure}[t]
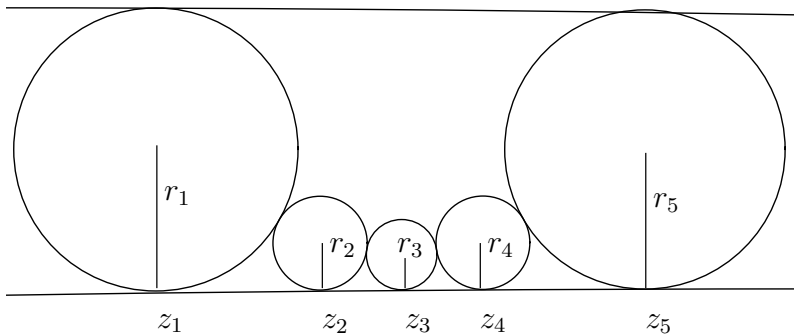

  \begin{center}
    \input aflower2.tex
    \caption{A flower with one central touching point at infinity.}
    \label{f.lines-circles}
  \end{center}
\end{figure}
{\em Proof.}
Mapping the point $z_6$ to infinity by a M\"obius transformation one obtains two
parallel lines and five touching circles as in Fig.\ref{f.lines-circles}. An
elementary computation yields
\begin{equation}                        \label{radii}
z_{k+1}-z_k= 2\sqrt{r_{k+1} r_k}, \qquad k=1,\ldots,4,
\end{equation}
where $r_k$ are the radii of the corresponding circles. Together with $r_1=r_5$
and $(z_5-z_6)/(z_6-z_1)=-1$ this implies (\ref{m-ratio}).

On the other hand, given arbitrary $r_1>0$ and ordered $z_1,\ldots,z_6$ satisfying
(\ref{m-ratio}), after normalizing $z_6=\infty$ formula (\ref{radii}) provides us with the radii
of the touching circles as in Fig.\ref{f.lines-circles}. This proves the equivalence of (i)
and (ii).

To show the equivalence of (ii) and (iii), define the M\"obius transformation $M$ through
$M(z_1)=z_4, M(z_2)=z_5, M(z_3)=z_6$.  Consider $z_*=M(z_4)$. The invariance of the cross-ratios
$q(z_1,z_2,z_3,z_4)$ $= q(z_4,z_5,z_6,z_*)$ implies the equivalence of (\ref{m-ratio}) and $z_*=z_1$.
The same proof holds for $M(z_5)=z_2$ and $M(z_6)=z_1$.
\vspace{5mm}

To each center touching point $z_k$ of a flower, one can associate a circle $S_k$ passing through 4
touching points $z_{k-1}, z_{k+1}, w_k, w_{k-1}$ of the flower containing $z_k$
(see Fig.\ref{f.s-circles}). Here $w_k$ is the touching point of
petals\footnote{The petals are labeled by the
corresponding touching points $z_k$.} ${\cal P}_{k+1}$ and ${\cal P}_k$.
Indeed, mapping the point $z_k$ by a M\"obius transformation to $\infty$, it is easy to see that
the points $z_{k-1}, z_{k+1}, w_k, w_{k-1}$ are mapped to vertices of a rectangle, thus lie on a
circle.  We call these circles {\em s-circles} of a flower.
\begin{figure}[ht]
  \begin{center}
    \online{%
    \centerline{\applet{cp.Flowerapplet.class}{280bp}{280bp}}
    \caption{A conformally symmetric flower. The green circles are the
      s-circles. If they all intersect in one point $p$ this one is
      shown red. Again you may drag the magenta-colored points and
      right mouse click over the circles gives you more viewer
      options. For help on the viewer click \protect\hyperlink{applethelp}{here}.}
    }{%
    \input afullFlower.tex
    \caption{A conformally symmetric flower.}
    }%
    \label{f.s-circles}
  \end{center}
\end{figure}

\begin{theorem}                             \label{t.F}
There exist a one-parameter family of flowers with the same center touching points. Moreover,
there exists a unique flower $F$ in this family, which satisfies the following equivalent
conditions:
\begin{itemize}
\item[(i)] $F$ is invariant with respect to a M\"obius involution $M$ with a fixed point $P$,
\item[(ii)] All s-circles of $F$ intersect in one point $P$.
\end{itemize}
\end{theorem}

We call the flower $F$ of the theorem {\em conformally symmetric}.

\online{In Figure~\ref{f.s-circles} one can view the whole family of flowers.}{
One can view the whole family of flowers at \cite{applets}.
}

\begin{figure}[ht]
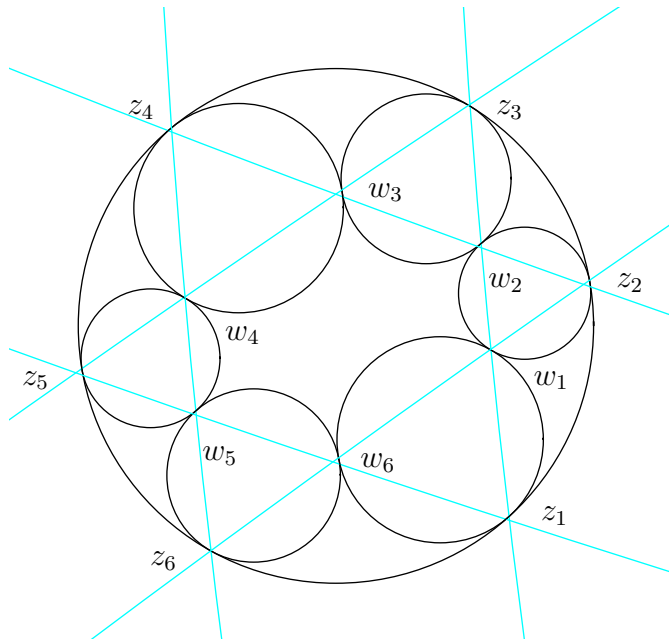

  \begin{center}
    \input anormalized.tex
    \caption{A normalized conformally symmetric flower.}            \label{f.star}
  \end{center}
\end{figure}

{\em Proof.}
Keeping the points $z_1,\ldots,z_5$ in Fig.\ref{f.lines-circles} fixed and varying $r_1$ one obtains
a one parameter family of flowers with the same touching central
points. Let us now
construct the flower $F$. The M\"obius involution of Proposition \ref{p.points-geometry}
preserves the central circle $C$. Consider the circles $C_k, \ i=1,2,3$, orthogonal to $C$
and passing through the pairs of points $\{z_k, z_{k+3}\}$. All these three circles intersect in
2 points $P$ and $P'$, which are the fixed points of $M$ lying inside and outside $C$,
respectively. By a M\"obius transformation, map the point $P$ to infinity.
The M\"obius involution $M$ becomes $M(z)=-z$ and the circles $C_1, C_2, C_3$ become straight lines
intersecting in the center of $C$. To construct the flower $F$, connect the $z_k$-points with even
(respectively, with odd)
labels by straight lines and consider their intersection points $w_k$
(see Fig.\ref{f.star}). The circles $C_k$ passing through the triples
$w_k, w_{k-1}, z_k$ touch at the points $w_k$. Let us prove this fact for $C_1$ and $C_2$. Indeed,
the triangles $\Delta(w_1, w_6, z_1)$ and $\Delta(z_3, z_5, z_1)$ are similar, therefore the
tangent lines to the circle $C_1$ at $w_1$ and to the circle $C$ at $z_3$ are parallel. The tangent
lines to $C_2$ at $w_1$ and to $C$ at $z_6$ are also parallel.
Since the points $z_3$ and $z_6$ are opposite on $C$, the circles $C_1$ and $C_2$ touch at $w_1$.
The circles $C_k$ are the petals of the desired flower $F$, which is obviously $M$-symmetric.
The s-circles of this flower are the straight lines $(z_k, z_{k+2})$. The latter obviously
intersect at infinity, thus all the s-circles of $F$ intersect in the fixed point $P$ of $M$.

The proof of $(ii)\Rightarrow (i)$ is similar. After mapping the point $P$ to infinity the s-circles
become straight lines and the flower is as in Fig.\ref{f.star}. Since
the circles in this figure touch, their
tangent lines at the points $z_k, z_{k+3}$ and $w_{k+1}$ are parallel. This implies that
 $z_k$ and $z_{k+3}$ are opposite points on $C$, and the flower is symmetric with respect to the
 $\pi$-rotation of $C$.
\vspace{5mm}

\begin{definition}
A hexagonal circle packing is called {\em conformally symmetric} or an
s-circle packing if it
consists of conformally symmetric flowers, i.e. the s-circles of each of its flowers intersect
in one point.
\end{definition}

\section{Analytic description of conformally symmetric circle packings}     \label{s.Analytic}

In this section we describe all conformally symmetric circle packings using the conformal
description of circle packings proposed by He and Schramm \cite{HS}.

To each central touching point $z_k$ of a flower one associates the cross-ratio\footnote{Note that
our normalization of $s_k$ differs from the one in \cite{HS}.}
\begin{equation}                    \label{s_i}
s_k:=q(z_k, z_{k-1}, w_{k-1}, w_k)=\frac{(z_k-z_{k-1})(w_{k-1}-w_k)}{(z_{k-1}-w_{k-1})(w_k-z_k)}.
\end{equation}
Mapping $z_k$ to $\infty$, one observes that tree other points in (\ref{s_i}) are mapped to
vertices of a rectangle, which implies that $s_k$ is purely imaginary. Moreover, the cross-ratios of
an immersed oriented flower are positive imaginary, $ -i s_k >0$.
Also note that
\begin{equation}                    \label{s-symmetry}
s_k= -q(z_{k+1}, z_{k-1}, w_{k-1}, z_k)= q(z_k, w_k, z_{k+1}, z_{k-1}),
\end{equation}
and that $s_k^2=q(z_{k+1}, z_{k-1}, w_{k-1}, w_k)$ is the
cross-ratio of the four touching points lying on the s-circle $S_k$.

\begin{lemma}
The cross-ratios $s_k$ of a flower satisfy the He-Schramm equation \cite{HS}
\begin{equation}                    \label{HSequation}
s_k +s_{k+2}+ s_{k+4}+ s_k s_{k+1} s_{k+2}=0
\end{equation}
for all $k\ {\rm mod}\ 6$.
\end{lemma}
{\em Proof.}
Let $m_k$ be the M\"obius transformation that takes $z_k, z_{k-1}, w_{k-1}$ to the points
$\infty, 0, 1$, respectively. By the definition of $s_k$ we have
\begin{eqnarray*}
s_k= q(z_k, z_{k-1}, w_{k-1}, w_k)=q(\infty, 0, 1, m_k(w_k)),\\
- s_k= q(z_{k+1}, z_{k-1}, w_{k-1}, z_k)= q(m_k(z_{k+1}), 0, 1, \infty),
\end{eqnarray*}
thus
$$
m_k(w_k)=1-s_k, \qquad m_k(z_{k+1})= -s_k.
$$
For $M_k:= m_{k+1}  m_1^{-1}$ this yields $M_k(-s_k)=\infty, M_k(\infty)=0, M_k(1-s_k)=1$ and,
finally,
$$
M_k=\left(
\begin{array}{cc}
0 & 1 \\
1 & s_k
\end{array}
\right),
$$
where the usual matrix notation for the M\"obius transformation is used. The equality of the
corresponding M\"obius transformations implies $M_3 M_2 M_1= \pm M_4^{-1}M_5^{-1}M_6^{-1}$,
which is
$$
\left(
\begin{array}{cc}
s_2 & 1 +s_1 s_2\\
1 +s_2 s_3 & s_1+s_3+s_1 s_2 s_3
\end{array}
\right)=\pm
\left(
\begin{array}{cc}
-s_4 -s_6 -s_4 s_5 s_6 & 1+ s_4 s_5\\
1+s_5 s_6 & - s_5
\end{array}
\right).
$$
Since the set of immersed flowers is connected and $s$'s do not vanish the sign in this equation is
the same for all flowers.
Taking all the circles with the same radius one checks that the correct sign is plus, which implies
the claim.
\vspace{5mm}

Given a hexagonal circle pattern it is convenient to associate its touching points as well as the
cross-ratios $s_k$ to the edges of the honeycomb lattice. Equation (\ref{HSequation}) is a partial
difference equation on the honeycomb lattice. The cross-ratios on the edges of each
hexagon satisfy (\ref{HSequation}). Moreover, it is easy to check that
the He-Schramm equation
is sufficient to guarantee the existence of the corresponding circle packing.

\begin{proposition}
Given a
positive-imaginary function $s:E\to i{\mathbb R}_+$ on the edges $E$ of the honeycomb lattice
satisfying (\ref{HSequation})
on each honeycomb, there exists unique (up to M\"obius transformation) immersed hexagonal circle
packing with the cross-ratios given by the corresponding values of $s$.
\end{proposition}

\begin{theorem}                     \label{t.s=s}
A circle flower is conformally symmetric if and only if its opposite cross-ratios $s_k$ are equal
\begin{equation}                    \label{constraint}
s_k=s_{k+3}\qquad\qquad (k\ {\rm mod}\ 6).
\end{equation}
\end{theorem}
{\em Proof.}
The property (\ref{constraint}) for conformally symmetric flowers follows from (i) of Theorem
\ref{t.F}. A simple computation with the flowers in Fig.\ref{f.star} shows that
the map $(s_1, s_2)$ of immersed conformally symmetric flowers to $(i{\mathbb R}_+)^2\ni (s_1, s_2)$
is surjective. Since a flower is determined through the $s$'s, the converse statement follows.
\vspace{5mm}

The general solution of (\ref{HSequation}, \ref{constraint}) on the whole $HL$ depends on three
arbitrary constants and can be given explicitly.
\online{
in Figure~\ref{f.s-circles} you can try how the circle packing changes
if one changes the three initial parameters.  
}{There is a JAVA applet that lets you explore this three parameter
family of circle packings interactively at \cite{applets}.
}

\begin{figure}[ht]
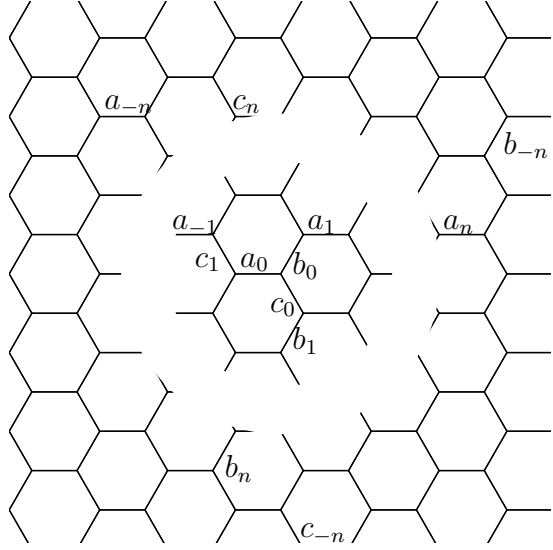

  \begin{center}
    \input agridLabel.tex
    \caption{Cross-ratios of conformally symmetric circle patterns.}
    \label{f.abc-labels}
  \end{center}
\end{figure}
\begin{theorem}
The general solution of (\ref{HSequation}, \ref{constraint}) is given by
\begin{eqnarray}
a_n & = & i\tan (\Delta n+\alpha),\nonumber\\
b_n & = & i\tan (\Delta n+\beta),               \label{solution}\\
c_n & = & i\tan (\Delta n+\gamma),\nonumber
\end{eqnarray}
where
$\Delta = -\alpha -\beta -\gamma$ and the cross-ratios $s_k$ on the edges of the hexagonal
lattice are labeled by $a_n, b_n, c_n$ as shown in Fig.\ref{f.abc-labels}.
\end{theorem}
{\em Proof.}
We start with a simple proof of the consistency of the following continuation of a solution of
(\ref{HSequation}, \ref{constraint}).
\begin{figure}[ht]
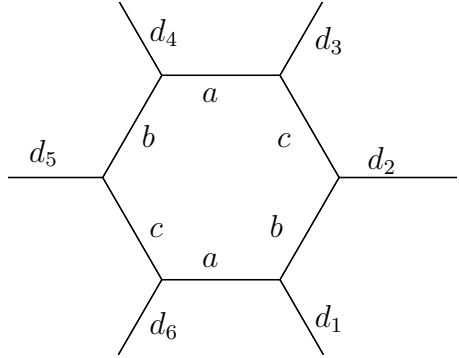

  \begin{center}
    \input ahoneyComb.tex
    \caption{Continuation of conformally symmetric $s$ about a honeycomb.}%
    \label{f.hired-honeycomb}
  \end{center}
\end{figure}
Given $s$ satisfying (\ref{HSequation}, \ref{constraint}) on a honeycomb $H$, i.e. $a,b,c$
in Fig.\ref{f.hired-honeycomb} satisfying
\begin{equation}
a+b+c+abc=0,                        \label{abc}
\end{equation}
and a value of $s$ on one of the edges attached to the honeycomb (
for example, $d_1$ in Fig.\ref{f.hired-honeycomb}), it can be uniquely extended to the full six
honeycombs $H_1,\ldots, H_6$ neighboring $H$. Indeed, (\ref{HSequation}, \ref{constraint}) yield
$$
b+d_1+d_2+b d_1 d_2=0,
$$
thus $d_2=M_1(d_1)$ is a M\"obius transformation of $d_1$. Passing once around the honeycomb
$H$ in this way one can check that (\ref{abc}) implies the monodromy M\"obius transformation
$M=M_6 \ldots M_1$ is the identity, thus this continuation implies no constraints on $d_1$.

Proceeding this way, one reconstructs $s$ on the whole lattice $HL$ from its values on
three adjacent edges ($a,b,d_1$ above). Then (\ref{HSequation}, \ref{constraint}) imply
$$
a_n+b_{-n}+c_1+a_n b_{-n} c_1=0,\qquad a_{n+1}+b_{-n}+c_0+a_{n+1} b_{-n} c_0=0
$$
and similar relations for other $a_n, b_n, c_n$. These identities become just
the addition theorem for the tangent function, implying the formulas
in (\ref{solution}), which can be
checked directly.
\online{
\begin{figure}[htbp]
  \begin{center}
    \centerline{\applet{cp.CSCPapplet.class}{360bp}{360bp}}
    \caption{Explore the conformally symmetric circle patterns. The
      evolve and shrink buttons increase and decrease the number of
      circles shown. The reset params button resets the sliders
      below. The three sliders control the three initial
      paprameters. They default to $\arctan(\sqrt{3})$ which corresponds
      to the standard packing. Right mouse click over the circle
      packing gives you a menu with viewer options. You might want to
      zoom out (using ``scale'' in the menu) at first. For help on the viewer click
      \protect\hyperlink{applethelp}{here}.} 
    
  \end{center}
\end{figure}
}{}
\section{The Doyle spirals}                     \label{s.Doyle}

\begin{figure}[ht]
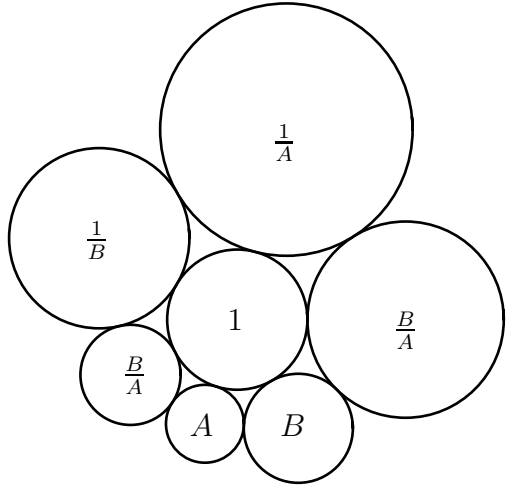

  \begin{center}
    \input aradii.tex
    \caption{Radii of a Doyle spiral with the normalized central radius $R=1$.}
    \label{f.DoyleRadii}
  \end{center}
\end{figure}

Denote by $R$ the radius of the center circle of a flower and by $R_k,\ k=1,\ldots,6$, the radii
of its petals. The Doyle spirals are characterized through the constraint (see \cite{BDS,CR} for a
complete analysis of Doyle spirals)
\begin{equation}                \label{RDoyle}
R_k R_{k+3}=R^2,\qquad R_k R_{k+2} R_{k+4}=R^3
\end{equation}
on the radii of the circles (see Fig.\ref{f.DoyleRadii} where the central radius is normalized to
be $R=1$). The Doyle spirals have two degrees of freedom (for example the ratios $R_1/R$ and $R_2/R$,
which are the same for the whole spiral) up to similarities.
\online{
  Figure~\ref{f.doyle.appl} You can experiment with the two radii and
  see how the corresponding Doyle spiral changes. 
}{Again, you can experiment with the two radii in a JAVA applet \cite{applets}.
}

\begin{proposition}
The Doyle spirals are conformally symmetric.
\end{proposition}
{\em Proof.}
The configurations of four touching circles with the radii $R, R_{k-1}, R_k, R_{k+1}$ and
with the radii $R_{k+3}, R_{k+4}, R, R_{k+2}$ differ by scaling.
This implies $s_k=s_{k+3}$ (use both (\ref{s_i}) and the second representation of $s_k$ in
(\ref{s-symmetry})) and the claim follows by Theorem \ref{t.s=s}.
\vspace{5mm}

\begin{theorem}                     \label{t.s=const}
The Doyle spirals and their M\"obius transformations can be characterized by the following two
equivalent properties:
\begin{itemize}
\item[(i)] The circle packing is conformally symmetric, and the corresponding solution of
(\ref{HSequation}) is ``constant''. It is of the
form (\ref{solution}) with $\alpha, \beta, \gamma \in (0, \frac{\pi}{2}), \
\alpha +\beta +\gamma=0 \ ({\rm mod} \ \pi)$ or, equivalently,
$$
a_n=a_0,\quad b_n=b_0,\quad c_n=c_0,\quad a_0, b_0, c_0\in i{\mathbb R}_+,\quad
 a_0+b_0+c_0+a_0 b_0 c_0=0.
$$
\item[(ii)] The whole circle packing is invariant with respect to the M\"obius involution
of each of its flowers.
\end{itemize}
\end{theorem}
{\em Proof.}
All the flowers of a Doyle spiral differ by scaling, which implies (i). Consider the Doyle spiral as
in Fig.\ref{f.DoyleRadii}. Computing the cross-ratios through the radii, one shows that the map
\begin{eqnarray*}
\{ (A,B)\in {\mathbb R}_+^2 \}  \to \{(a,b,c)\in (i{\mathbb R}_+)^3 : a+b+c+abc=0\}.
\end{eqnarray*}
is surjective, thus (i) characterizes the Doyle spirals and their M\"obius transforms.
The proof of the equivalence $(i)\Leftrightarrow(ii)$ is elementary and is left to the reader.
\vspace{5mm}

It is an open problem whether the Doyle spirals are the only entire circle packings. Formulas
(\ref{solution}) imply that it is possible to have all cross-ratios
being positive imaginary (necessary condition
for entireness) only when $\Delta=0$.
\begin{corollary}
Doyle spirals are the only entire conformally symmetric circle patterns.
\end{corollary}
\online{
\begin{figure}[htbp]
  \begin{center}
    \centerline{\applet{cp.DoyleApplet.class}{360bp}{360bp}}
    \caption{Explore the Doyle spirals: The sliders change the radii
      equally labeled circles. Right mouse click over the circle
      packing gives you a menu with viewer options. For more help
      click \protect\hyperlink{applethelp}{here}}
    \label{f.doyle.appl}
  \end{center}
\end{figure}
}{}
\section{Airy functions as continuous limit}

Because of the property (\ref{RDoyle}), Doyle spirals are interpreted as a discrete exponential
function.

In the conformal setting this interpretation can also be easily observed. Indeed, let $P^\epsilon$ be
a family of circle packings approximating a holomorphic mapping in the limit $\epsilon \to 0$. In
\cite{HS} He and Schramm investigated the behavior of the cross-ratios $s_k$ in this limit:
$$
s_k = i\sqrt{3}(1+ \epsilon^2 h_k^\epsilon),
$$
where $h_k$ is called the discrete Schwarzian derivative (Schwarzian) of $P^\epsilon$ at the
corresponding edge of the hexagonal lattice. The discrete Schwarzians converge to the
Schwarzian derivative
\begin{equation}                            \label{S(f)}
S(f):=\left(\frac{f''}{f'}\right)' -\frac{1}{2}\left(\frac{f''}{f'}\right)^2
\end{equation}
of the corresponding holomorphic mapping.
More precisely, there exist continuous limits
$$
a=\lim_{\epsilon\to 0} h_1^\epsilon,\quad
b=\lim_{\epsilon\to 0} h_2^\epsilon,\quad
c=\lim_{\epsilon\to 0} h_3^\epsilon
$$
for the smooth functions\footnote{Note that $\lim_{\epsilon\to 0} h_k^\epsilon=
\lim_{\epsilon\to 0} h_{k+3}^\epsilon$} $a, b, c$. Because of (\ref{HSequation}) these functions
satisfy
\begin{equation}                            \label{a+b+c=0}
a+b+c=0
\end{equation}
at each point. The Schwarzian equals
$$
S(f)=4(a+q^2 b+ q c),\qquad q=e^{2\pi i/3},
$$
and, using (\ref{a+b+c=0}), this also yields
\begin{equation}                            \label{a=Re(f)}
6a={\rm Re} (S(f)),\quad 6b={\rm Re} (q S(f)),\quad 6c={\rm Re} (q^2 S(f)).
\end{equation}
We see that, due to Theorem \ref{t.s=const}, Doyle spirals correspond to holomorphic functions with
constant Schwarzian derivative $S(f)={\rm const}$. The general solution of the last equation is the
exponential function and its M\"obius transformations.

It is natural to ask which holomorphic functions correspond to general
conformally symmetric circle packings. In Fig.\ref{f.abc-labels} one observes that each of the
cross-ratios $a_n, b_n, c_n$ is constant along one lattice direction. For the functions $a,b,c$ above,
this implies
\begin{equation}                            \label{a=a(x)}
a=a({\rm Re} (z)),\qquad b=b({\rm Re} (q z)), \qquad c=c({\rm Re} (q^2 z)),
\end{equation}
where $z$ is the complex coordinate. Comparing (\ref{a=Re(f)}) and (\ref{a=a(x)}) we see that the Schwarzian
is a linear function of $z$:
\begin{equation}                            \label{S-linear}
S(f)= Az+B, \qquad A\in {\mathbb R},\ B\in {\mathbb C}.
\end{equation}
Equation (\ref{S-linear}) can be easily solved by standard methods.
The general solution of
$
S(f)=u(z)
$
with holomorphic $u(z)$ is given by
$
f(z):=\psi_1/\psi_2,
$
where $\psi_1(z)$ and $\psi_2(z)$ are two independent solutions of the linear differential equation
$
\psi''= u(z) \psi.
$

\begin{figure}[ht]
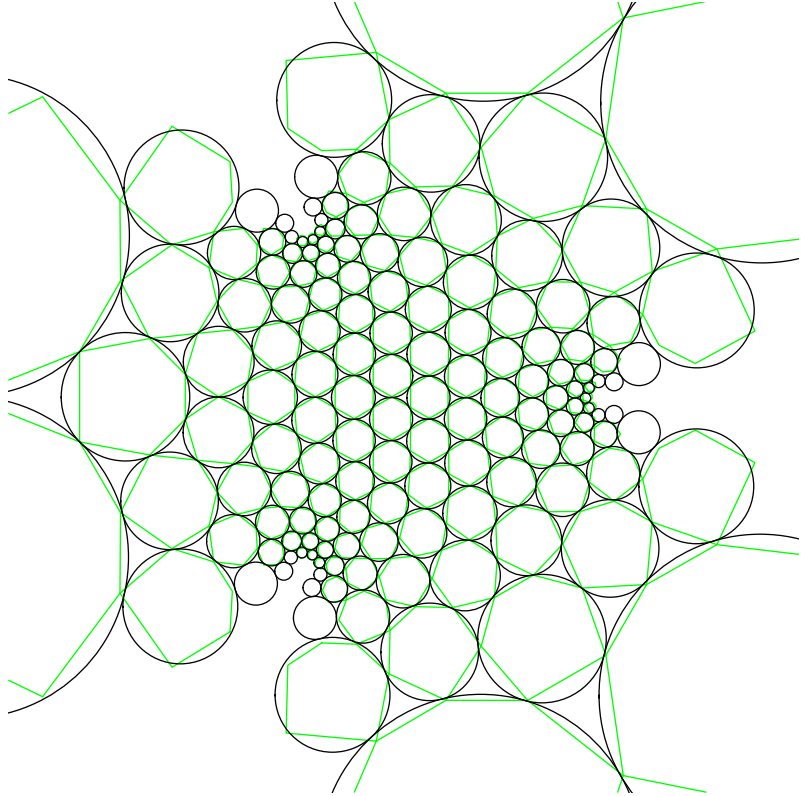

  \leavevmode
  \begin{center}
    \input aAiryBoth.tex
    \caption{A conformally symmetric circle packing (with
      $\alpha=\beta=\gamma$ in (\ref{solution})) and its smooth
      counterpart. The vertices of the
      hexagons are the images of the points of a standard hexagonal grid under the map
      $f$ from (\ref{e.fDef}).}
    \label{f.Airy}
  \end{center}
\end{figure}

By a shift and scaling of the variable $z$, equation (\ref{S-linear}) with $A\not= 0$ can be brought to
the form
\begin{equation}                            \label{S=z}
S(f)=z.
\end{equation}
Solutions of the corresponding linear equation
\begin{equation}                            \label{Airy}
\psi''=z\psi
\end{equation}
are the Airy functions ${\rm Ai}(z)$ and ${\rm Bi}(z)$. On the real line the first one is given by
\cite{SO}
$$
{\rm Ai}(x)=\frac{1}{\pi}\int_0^\infty \cos (xt+\frac{t^3}{3})dt,
$$
and the second one is related to it by
$$
{\rm Bi}(z)=iq^2 {\rm Ai}(q^2 z)-iq {\rm Ai}(q z).
$$
In the corresponding M\"obius class of solutions of (\ref{S=z}) it is natural to choose
\begin{equation}
f(z):= \frac{{\rm Bi}(z)-\sqrt{3}{\rm Ai}(z)}{{\rm Bi}(z)+\sqrt{3}{\rm Ai}(z)},
\label{e.fDef}
\end{equation}
which is the most symmetric one,
$
f(qz)=qf(z).
$
The corresponding circle packing, symmetric with respect to the rotation $z\to qz$,
is shown in Fig.\ref{f.Airy}.

\online{\newpage
\hypertarget{applethelp}{\section*{Using the circle packing viewer}}

Pressing the right mouse button (or holding $<META>$-key while pressing
the mouse on one button systems) over the viewers main window will
bring up a pop up menu with the following items
\begin{description}
\item[translate] After selecting this you may drag the whole viewed
  configuration
\item[scale]  After selecting this dragging the mouse up or down will
  zoom in or out
\item[rotate] After selecting this dragging the mouse up or down will
  rotate the image
\item[moebiustransform] this option allows to move ``0'', ``1'' and
  ``$\infty$'' freely and thus apply an arbitrary M\"obius
  transformation. These three points are shown with label (two in blue
  and one in red). Mouse dragging will drag the red one and a single
  click will step through the three points selecting the next one. You
  may also directly pick and drag any of them. If one holds down
  $<Control>$-key while dragging instead of the selected point its
  inverse will be moved. This way one can get the point ``$\infty$'' away
  from $\infty$ in the beginning: Initially only ``0'' and ``1'' are
  visible. Click until both are blue. Then $<Control>$-drag until the
  red ``inf'' labeled point moves into the visible area. Then you can
  release $<Control>$ and drag ``inf'' directly.
\item[invert on circle] lets you invert the whole configuration on an
  arbitrary circle. Click-dragging will open a yellow
  circle. Releasing the mouse button will do the inversion.
\item[save view] lets you store the current view on the configuration
  for later return.
\item[restore view] lets you successively restore any previously saved
  views.
\item[reset view] resets the view to the initial settings.
\item[move points] in applets where points may be moved around this
  option makes the points drag-able again when one has switched to
  other options in between.
\end{description}
}{
\section{About the interactive version of this paper}\label{s.online}
Since the families of circle packings discussed in this paper have a
finite (and even small) number of parameters it seemed to be natural to
look for a way to visualize the whole families.

Except the Doyle spirals the families are only defined modulo
an arbitrary M\"obius transformation. Therefore it should be possible to view
the packings with a M\"obius transformation applied, so one could look
at them through ``M\"obius glasses''.

The outcome of this is an interactive
version of this paper. This version includes 
some java applets that let you experiment with the circle packings
directly. In particular there are applets to illustrate the families
of circle packing flowers, the whole class of conformally symmetric
circle packings, and the special case of Doyle spirals.

This interactive version renders a dvi file
inside a java applet. It needs a web browser that includes a java
vm. Since flipping pages might get slow on old machines and java vm's that have no just in time
compiler, we also provide a page that shows the applets
that are missing in this version only. You will find the paper and the
applets at
\\[0.3cm]
\centerline{\href{http://www-sfb288.math.tu-berlin.de/Publications/online/}%
{\bf http://www-sfb288.math.tu-berlin.de/Publications/online/}}
\\[0.3cm]
}

\end{document}